\documentclass[11pt]{amsart}
\usepackage{amsmath}

\theoremstyle{plain}
\newtheorem{theorem}{Theorem}[section]
\newtheorem{proposition}[theorem]{Proposition}
\newtheorem{lemma}[theorem]{Lemma}
\newtheorem{corollary}[theorem]{Corollary}
\newtheorem{fact}[theorem]{Fact}
\newtheorem*{claim}{Claim}
\newtheorem*{question}{Question}

\newcommand{\Zn}{{\mathbb Z}(n)}
\newcommand{\Zpn}{{\mathbb Z}(p^n)}

\newcommand{\Zpks}{{\mathbb Z}(p^{k+1})}
\newcommand{\Zpi}{{\mathbb Z}(p^\infty)}
\newcommand{\Zp}{{\mathbb Z}_{(p)}}
\newcommand{\bQ}{\mathbb Q}

\DeclareMathOperator{\End}{End}

\title
[Theory of square-like abelian groups is decidable]
{Theory of square-like abelian groups\\ is decidable}

\author{Oleg Belegradek}
\address{Department of Mathematics\\
Istanbul Bilgi University\\
80370 Dolapdere--Istanbul, Turkey}
\email{\tt olegb@bilgi.edu.tr}


\date{February 9, 2006}

\begin{document}

\begin{abstract}
A group is called square-like if it is universally
equivalent to its direct square.
It is known that the class of all square-like groups
admits an explicit first order axiomatization
but its theory is undecidable. We prove that the theory
of square-like \emph{abelian} groups is decidable.
This answers a question posed by D. Spellman.

\end{abstract}

\maketitle

\section*{Introduction}
A group $G$ is called \emph{discriminating}~\cite{BMR}
if every group separated by $G$ is discriminated by $G$.
Here $G$ is said to separate (discriminate)
a group $H$ if for any non-identity element
(finite set of non-identity elements) of $H$
there is a homomorphism from $H$ to $G$
which does not map the element (any element of the set) to the identity.
A group $G$ is discriminating iff
$G$ discriminates $G^2$~\cite{BMR}. In particular,
if $G$ embeds $G^2$ then
$G$ is discriminating.

A group $G$ is called \emph{square-like}~\cite{FGMS2}
if the groups $G^2$ and $G$ are universally equivalent.
Any discriminating group is square-like~\cite{FGMS1}.
The notions of discriminating and square-like group were studied
in \cite{BMR,B2,FGMS1,FGMS2,FGS1,FGS2,FGS3,FGS4}.

The class of square-like groups
is first order axiomatizable~\cite{FGMS2}, and
the theory of the class is computably enumerable;
an explicit first order axiom system was
suggested in \cite{B1,B2}, and also presented in~\cite{FGS3}.
In \cite{FGMS2} square-like abelian groups
were characterized in terms of Szmielew invariants.

The subclass of discriminating groups is not first order
axiomatizable~\cite{FGMS2}. Every square-like group is
elementarily equivalent to a discriminating group \cite{B2,FGS2};
so the class of square-like groups is the axiomatic closure
of the class of discriminating groups.

The theory of square-like groups is undecidable~\cite{B2, FGS2}.
The argument in~\cite{FGS2}
is based on the obvious observation that
any group embeds in a discriminating group,
and so the universal theory of square-like groups
coincide with the universal theory of all groups.
The latter is undecidable because there exist
finitely presented groups with unsolvable word problem.
In~\cite{B2} a discriminating group that interprets
the ring of integers is constructed;
any theory that has the group as a model
(and, in particular, the theory of square-like groups)
is undecidable.

The main result of the present paper is that
the theory of square-like abelian groups is decidable.
This answers a question posed by Dennis Spellman~\cite{S}.
As a byproduct, we found characterizations
of discriminating and square-like Szmielew groups.

\section{Preliminaries}
Here we collect some known definitions and facts we will use in the proofs.

\begin{fact}\cite[Proposition 1]{BMR}\label{square}
A group $G$ is discriminating iff $G$ discriminates $G^2$.
In particular,
$G$ is discriminating if $G$ embeds $G^2$.
\end{fact}

\begin{fact}\cite[Proposition 2]{BMR}\label{direct-sum}
The direct product \textup{(}restricted or not\textup{)}
of any family of discriminating groups is
a discriminating group.
\end{fact}

\begin{fact}\cite[Proposition 3]{BMR}\label{torsion-free}
Any torsion-free abelian  group is discriminating.
\end{fact}

\begin{fact}\cite[Lemma 2.1]{FGMS1}
\label{discriminating-is-square-like}
Any discriminating group is square-like.
\end{fact}

\begin{fact}\cite[Theorem 3]{FGMS2}\label{square-like-axioms}
The class of square-like groups is first order
axiom\-atizable.
\end{fact}

\begin{fact}\cite[Proposition 3.5]{B2}\label{end-invariant}
Any $\End(G)$-invariant subgroup of
a discriminating group $G$ is
trivial or infinite.
\end{fact}

Let $A$ be an abelian group. For a positive integer $n$
we denote
$$nA=\{na: a\in A\},\quad A[n]=\{a\in A: na=0\},$$
and write $\delta(A)$ for the largest divisible subgroup
of $A$. We write $nA[k]$ for $(nA)[k]$.
The subgroups $nA$, $A[n]$, $nA[k]$, and $\delta(A)$ are
$\End(A)$-invariant.
We write $A^{(\kappa)}$ for the direct sum of $\kappa$
copies of $A$.

We write $\bQ$ for the additive group of all rational numbers,
and $\Zp$ for the additive group of rational numbers with denominator
not divisible by a prime $p$.
We write $\Zn$
for the cyclic group of order $n$,
and $\Zpi$ for the Pr\"{u}fer $p$-group.

A \emph{Szmielew group} is defined to be an abelian group
of the form
\begin{equation}
\bigoplus_{p\ \text{prime}}[
\bigoplus_{n>0}
\Zpn^{(\kappa_{p,n-1})}\oplus
\Zpi^{(\lambda_p)}\oplus
\Zp^{(\mu_p)}]
\oplus\bQ^{(\nu)}\tag{$\star$}
\end{equation}
where
$\kappa_{p,n-1}$,
$\lambda_p$,
$\mu_p$, $\nu$ are cardinals $\le\omega$.

For a prime $p$, we call a Szmielew group of the form
\begin{equation*}
\bigoplus_{n>0}
\Zpn^{(\kappa_{p,n-1})}\oplus
\Zpi^{(\lambda_p)}\oplus
\Zp^{(\mu_p)}
\oplus\bQ^{(\nu)}
\end{equation*}
a \emph{$p$-Szmielew} group.

\begin{fact}\cite[Lemma A.2.3]{H}\label{equivalent-to-Szmielew}
Every abelian group is elementarily equivalent to a
Szmielew group.
\end{fact}

Let $p$ be a prime,
and $n,k<\omega$.
Let
$\Phi_k(p,n)$ and $\Phi^k(p,n)$
be the sentences that say about an abelian group $B$ that
$$
\dim_p(p^nB[p]/p^{n+1}B[p])=k\quad\text{and}\quad
\dim_p(p^nB[p]/p^{n+1}B[p])>k,
$$
$\Theta_k(p,n)$ and $\Theta^k(p,n)$
be the sentences that say that
$$
\dim_p(p^nB[p])=k\quad\text{and}\quad
\dim_p(p^nB[p])>k,
$$
$\Gamma_k(p,n)$ and $\Gamma^k(p,n)$
be the sentences that say that
$$
\dim_p(p^nB/p^{n+1}B)=k\quad\text{and}\quad
\dim_p(p^nB/p^{n+1}B)>k,
$$
$\Delta_k(p,n)$ and $\Delta^k(p,n)$
be the sentences that say that
$$
|p^nB|=k
\quad\text{and}\quad
|p^nB|>k.
$$
The sentences defined above are called the
Szmielew invariant sentences.
Note that
$|B|=k$ and
$|B|>k$ can be expressed as
$\Delta_k(p,0)$ and $\Delta^k(p,0)$,
for any prime $p$.

\begin{fact}\cite[Section A.2]{H}\label{Szmielew-invariants}
If $A$ is the Szmielew group $(\star)$ then
\begin{itemize}
\item
$A\models\Phi_k(p,n)$\quad iff\quad $\kappa_{p,n}=k$,
\smallskip

\item
$A\models\Phi^k(p,n)$\quad iff\quad $\kappa_{p,n}>k$,
\smallskip

\item
$A\models\Theta_k(p,n)$\quad iff\quad
$\lambda_p+\kappa_{p,n}+\kappa_{p,n+1}+\dots=k$,
\smallskip

\item
$A\models\Theta^k(p,n)$\quad iff\quad
$\lambda_p+\kappa_{p,n}+\kappa_{p,n+1}+\dots>k$,
\smallskip

\item
$A\models\Gamma_k(p,n)$\quad iff\quad
$\mu_p+\kappa_{p,n}+\kappa_{p,n+1}+\dots=k$,
\smallskip

\item
$A\models\Gamma^k(p,n)$\quad iff\quad
$\mu_p+\kappa_{p,n}+\kappa_{p,n+1}+\dots>k$.
\end{itemize}
\end{fact}
\smallskip

\begin{fact}\cite[Theorem A.2.7]{H}\label{translation-Szmielew}
Every sentence of the first order language of abelian groups
is equivalent, modulo the theory of abelian groups,
to a positive Boolean combination of Szmielew invariant
sentences.
\end{fact}

\begin{fact}\cite[Theorem A.2.7]{H}\label{Szmielew-equivalence}
Two abelian groups are elementarily equivalent
iff they satisfy the same Szmielew invariant sentences.
\end{fact}

Abusing terminology,
we call a sentence of the language of abelian groups \emph{consistent}
if it is true in some abelian group.
By Fact~\ref{equivalent-to-Szmielew},
a sentence is consistent iff it holds in some Szmielew group.

\begin{fact}\cite[Theorem A.2.8]{H}\label{Szmielew-consistency}
There is an algorithm that, given
a finite conjunction of Szmielew invariant sentences,
decides whether it holds in some Szmielew group.
\end{fact}

Facts \ref{translation-Szmielew} and \ref{Szmielew-consistency}
are main ingredients of a proof of
the Szmielew theorem
on decidability of the theory of abelian groups;
actually, they immediately imply the result.
Indeed, given a sentence $\phi$,
by Fact \ref{translation-Szmielew}
and  computable enumerability of the theory of abelian groups,
we can effectively find
a positive Boolean combination $\theta$ of Szmielew invariant
sentences that is equivalent to $\neg\phi$, modulo the theory.
A sentence $\phi$ is not in the theory iff $\theta$
is consistent; the latter can be effectively checked, by
Fact~\ref{Szmielew-consistency}.

We will use a similar method in our proof of
decidability of the theory of square-like abelian groups.

\section{Discriminating and square-like Szmielew groups}

Let $A$ be the Szmielew group $(\star)$.
For a prime $p$, let $I_p=\{n:\kappa_{p,n-1}>0\}$.
In case when the set $I_p$ is finite and nonempty,
$l_p$ denotes its maximal element;
clearly,
$\kappa_{p,l_p-1}>0$.

\begin{proposition}\label{discr}
The following are equivalent:
\begin{enumerate}
\item
$A$ is discriminating;
\item
for any prime $p$ one of the following holds:
\begin{itemize}
\item[(i)]
$\lambda_p=\omega$,
\item[(ii)]
$\lambda_p=0$, and if $I_p$ is finite and nonempty then
$\kappa_{p,\,l_p-1}=\omega$.
\end{itemize}
\end{enumerate}
\end{proposition}

\begin{proof}
(1)$\Rightarrow$(2).
Suppose (1). Let $p$ be a prime.
The subgroup $\delta(A)\cap A[p]$
is ${\rm End}(A)$-invariant, and hence is trivial or infinite,
by Fact~\ref{end-invariant}.
Then $\lambda_p$ is 0 or $\omega$.
Suppose  $\lambda_p=0$, and $I_p$ is finite and nonempty.
Then the ${\rm End}(A)$-invariant subgroup
$p^{l_p-1}A[p]$
is nontrivial and hence infinite, again by
Fact~\ref{end-invariant}.
Then $\kappa_{p,\,l_p-1}=\omega$.

(2)$\Rightarrow$(1).
Suppose (2). Then for any prime $p$ the group
$$\bigoplus_{n>0}
\Zpn^{(\kappa_{p,n-1})}\oplus
\Zpi^{(\lambda_p)}$$
embeds it square.
So $A=B\oplus C$, where $B$ embeds $B^2$, and $C$ is torsion-free.
By Facts~\ref{square}, \ref{torsion-free}, and \ref{direct-sum},
$A$ is discriminating.
\end{proof}

\begin{proposition}\label{square-like}
The following are equivalent:
\begin{enumerate}
\item
$A$ is square-like;
\item
for any prime $p$ one of the following holds:
\begin{itemize}
\item[(i)]
$\lambda_p=\omega$,
\item[(ii)]
$\lambda_p=0$, and if $I_p$ is finite and nonempty then
$\kappa_{p,\,l_p-1}=\omega$,
\item[(iii)]
$0<\lambda_p<\omega$, and $I_p$ is infinite.
\end{itemize}
\end{enumerate}
\end{proposition}

\begin{proof}
$(1)\Rightarrow(2)$.
Suppose (2) fails. Then, for some
prime $p$, (i), (ii), (iii) all fail.
There are two possibilities:
\begin{enumerate}
\item[(a)]
$\lambda_p=0$, the set $I_p$ is finite, nonempty, and
$\kappa_{p,\,l_p-1}<\omega$,
\item[(b)]
$0<\lambda_p<\omega$, and the set $I_p$ is finite.
\end{enumerate}
Suppose (a). Let $\kappa=\kappa_{p,\,l_p-1}$.
We have
$$|p^{l_p-1}A[p]|=p^\kappa,\qquad
|p^{l_p-1}A^2[p]|=p^{2\kappa}.
$$
Suppose (b). Put $l=l_p$ if
$I_p\ne\emptyset$, and
$l=0$ otherwise. We have
$$|p^lA[p]|=p^{\lambda_p},\qquad
|p^lA^2[p]|=p^{2\lambda_p}.
$$
For any positive integers $s$ and $t$ there is an existential
sentence that says about an abelian group $B$ that
$|sB[p]|\ge t$.
Therefore in both cases (a) and (b) the groups
$A$ and $A^2$ are not universally equivalent,
and so (1) fails.

$(2)\Rightarrow(1)$.
Suppose (2). Let $A'$ be the Szmielew group obtained
from $A$ by replacing
$$\bigoplus_{n>0}
\Zpn^{(\kappa_{p,n-1})}\oplus
\Zpi^{(\lambda_p)}$$
with
$$\bigoplus_{n>0}
\Zpn^{(\kappa_{p,n-1})},$$
for all $p$ satisfying (3).
Then $A'$ is discriminating, by Proposition~\ref{discr}.
Hence $A'$ is square-like, by Fact~\ref{discriminating-is-square-like}.
It is easy to check that
$A$ and $A'$ satisfy the same Szmielew invariant sentences;
therefore, by Fact~\ref{Szmielew-equivalence},
$A\equiv A'$. Then, by Fact~\ref{square-like-axioms},
the group $A$ is square-like, too.
\end{proof}

\begin{corollary}\label{equivalent-to-discriminating}
Any square-like abelian group is elementarily equivalent to
a discriminating Szmielew group.
\end{corollary}

\begin{proof}
Let $B$ be a square-like abelian group.
By Fact~\ref{equivalent-to-Szmielew},
$B$ is elementarily equivalent to
a Szmielew group $A$.
By Fact~\ref{square-like-axioms},
$A$ is square-like.
The argument at the end of the proof of Proposition~\ref{square-like}
shows that $A$ is elementarily equivalent to a discriminating
Szmielew group $A'$.
\end{proof}

\section{Main result}

\begin{theorem}\label{decidable}
The theory of square-like abelian groups is decidable.
\end{theorem}

\begin{proof}
We need to find an algorithm which,
given a sentence $\phi$ of the language of abelian groups,
decides whether $\phi$ is true in some square-like abelian group,
or, equivalently by
Corollary~\ref{equivalent-to-discriminating},
in some discriminating Szmielew group.
By Fact~\ref{translation-Szmielew},
$\phi$ is equivalent, modulo the theory
of abelian groups, to a positive Boolean combination $\theta$
of Szmielew invariant sentences. Since the theory of abelian groups
is computably enumerable, $\theta$ can be found effectively.
We may assume that $\theta$ is $\bigvee_i\theta_i$, where
each $\theta_i$ is
a conjunction of
finitely many Szmielew invariant sentences.
So it suffices to prove

\begin{claim}
There exists an algorithm that,
given a consistent conjunction $\psi$
of finitely many Szmielew invariant sentences,
decides whether $\psi$ holds in some
discriminating Szmielew group.
\end{claim}

For a prime $p$, we call a conjunction
of formulas  of the forms
\begin{gather*}
\Phi_k(p,n),\
\Theta_k(p,n),\
\Gamma_k(p,n),
\Delta_k(p,n),\\
\Phi^k(p,n),\
\Theta^k(p,n),\
\Gamma^k(p,n),\
\Delta^k(p,n)
\end{gather*}
a $p$-conjunction.
To prove the Claim, we show that

\begin{enumerate}
\item[(A)]
\emph{there exists an algorithm that,
given a prime $p$ and a consistent
$p$\nobreakdash-conjunction $\psi$,
decides whether $\psi$ holds in some
discriminating $p$\nobreakdash-Szmielew group}, and
\smallskip

\item[(B)] \emph{the Claim follows from} (A).
\end{enumerate}

First we show (B): assuming (A), we prove the Claim.

Let $\psi$ be a conjunction of Szmielew invariant sentences,
which holds in a Szmielew group~$A$.
We have $\psi=\bigwedge_p\psi_p$, where $p$ runs over
a finite set of primes, and $\psi_p$ is a $p$-conjunction.
There are three possibilities:
\begin{enumerate}
\item[(a)]
$\psi$ has no conjuncts of the form $\Delta_k(p,n)$;
\item[(b)]
$\psi$ has some conjuncts
$\Delta_k(p,n)$ and
$\Delta_l(q,m)$ with $p\ne q$;
\item[(c)]
$\psi$ has a conjunct $\Delta_k(p,n)$, but has no
conjuncts $\Delta_l(q,m)$ with $p\ne q$.
\end{enumerate}

The following three lemmas prove (B).

\begin{lemma}\label{no-delta}

Assume \textup{(a)}. The following are equivalent:
\begin{enumerate}
\item[(i)]
$\psi$ holds in some
discriminating Szmielew group,
\item[(ii)]
for all $p$ the sentence
$\psi_p$ holds in some
discriminating $p$-Szmielew group.
\end{enumerate}
\end{lemma}

\begin{proof}
Suppose (i).
We have $A=\oplus_pA(p)$,
where $A(p)$ is a $p$-Szmielew group.
Let $p$ be a prime.
Then $A(p)\oplus\bQ$ is a discriminating $p$-Szmielew group,
by Proposition~\ref{discr}. Also,
$A(p)\oplus\bQ\models\psi_p$ because of (a).
So (ii) holds.

Suppose (ii). For every prime $p$ choose
a discriminating $p$-Szmielew group $A(p)$ in which
$\psi_p$ holds.
By Proposition~\ref{discr}, the Szmielew group
$A=\oplus_pA(p)$ is discriminating.
For every $p$
we have $A\models\psi_p$, because
$A(p)\models\psi_p$ and $\psi$ satisfies (a).
Therefore $A\models\psi$. So (i) holds.
\end{proof}

\begin{lemma}\label{exp}
Let $B$ be a discriminating abelian group.
\begin{enumerate}
\item
If $\Delta_k(p,n)$ or $\neg\Delta^k(p,n)$
holds in $B$ then $p^nB=0$.
\smallskip

\item
Assume \textup{(b)}. If $B\models\psi$ then $B=0$.
\end{enumerate}
\end{lemma}

\begin{proof}
(1) The subgroup $p^nB$ is
${\rm End}(B)$-invariant and finite of order at most~$k$.
By Fact~\ref{end-invariant},
the result follows.

(2) By (1), $p^nB=q^mB=0$, and hence $B=0$.
\end{proof}

Thus, for any $\psi$ with (b),
in order to decide whether
there is a discriminating Szmielew group that satisfies $\psi$,
we need to decide whether
$\psi$ holds in the trivial group,
which can be done effectively.

\begin{lemma}\label{delta}
Assume \textup{(c)}. Then
$\psi$ holds in some
discriminating Szmielew group if and only if
\begin{itemize}
\item[(i)]
For any $q\ne p$ and $l>0$, in
$\psi$ there are no conjuncts of the forms
$$\Phi^l(q,m),\ \Theta^l(q,m),\ \Gamma^l(q,m),\
\Phi_l(q,m),\
\Theta_l(q,m),\
\Gamma_l(q,m);$$
\item[(ii)]
For any $q\ne p$, in $\psi$ there are no conjuncts of the forms
$$\Phi^0(q,m),\
\Theta^0(q,m),\
\Gamma^0(q,m);$$
\item[(iii)]
the $p$-conjunction
$$\psi_p\wedge\bigwedge\{\Delta^s(p,0):s\in S\}$$
holds in some discriminating $p$-Szmielew group,
where $S$ is the set of all $s$ such that
$\Delta^s(q,m)$
is a conjunct of $\psi$, for some $q\ne p$
and some $m$.
\end{itemize}
\end{lemma}

\begin{proof}
First suppose that $\psi$ holds in a discriminating Szmielew group~$A$.
By (c) and Lemma~\ref{exp}\,(1),
$p^nA=0$, and so $A$ is a $p$-Szmielew group.
Therefore (i) and (ii) hold.
Let $s\in S$. Then for some $m$ and $q\ne p$ we have
$A\models\Delta^s(q,m)$, that is, $|q^mA|>s$.
As $p^nA=0$, we have $q^mA=A$; thus $|A|>s$.
Then $A\models\Delta^s(p,0)$. So (iii) holds.

Now suppose (i)--(iii) hold. By (iii) there is
a discriminating $p$\nobreakdash-Szmie\-lew group $A$ in which
$\psi_p$ and $\{\Delta^s(p,0):\ s\in S\}$ are true.
We show that $A\models\psi$.
Since
$\Delta_k(p,n)$ is a conjunct of $\psi$,
we have $p^nA=0$, by Lemma~\ref{exp}\,(1).
As $A$ is a $p$-Szmielew group, all the sentences
$\Phi_0(q,m)$,
$\Theta_0(q,m)$,
$\Gamma_0(q,m)$ with $q\ne p$ hold in $A$.
Due to (i) and (ii),
it remains to show that if
$\Delta^s(q,m)$
is a conjunct of $\psi$, where $q\ne p$,
then it holds in $A$.
Suppose not.
Then $q^mA=0$, by Lemma~\ref{exp}\,(1).
Therefore $A=0$, contrary to
$A\models\Delta^s(p,0)$.
\end{proof}

Now we prove (A).
From now on, \emph{let $p$ be a fixed prime,
and $\psi$ be a $p$\nobreakdash-conjunction which holds
in some Szmielew group $A$.}
We will show how to decide
whether $\psi$ holds in some discriminating
$p$\nobreakdash-Szmielew  group.

There are four possibilities:

\begin{enumerate}
\item[(a)]
$\psi$ has a conjunct
$\Delta_k(p,n)$ with $k\ne 1$;
\item[(b)]
$\psi$ has a conjunct
$\Theta_k(p,n)$ with $k>0$;
\item[(c)]
$\psi$ has no conjuncts of the forms
$\Delta_k(p,n)$ and
$\Theta_k(p,n)$;
\item[(d)]
$\psi$ has a conjunct $\Delta_1(p,n)$ or $\Theta_0(p,n)$,
but (a) and (b) fail.
\end{enumerate}

\begin{lemma}
If \textup{(a)} then
$\psi$ fails in  every discriminating abelian group.
\end{lemma}

\begin{proof}
Suppose
$\psi$ holds in an abelian group $B$.
Then $|p^nB|=k\ne 1$,
and so $p^nB$ is a nontrivial finite $\End(B)$-invariant
subgroup.
Therefore $B$ is not discriminating, by Fact~\ref{end-invariant}.
\end{proof}

\begin{lemma}
If \textup{(b)} then
$\psi$ fails in every discriminating Szmielew group.
\end{lemma}

\begin{proof}
Suppose $A\models\psi$, and $A$ is
a discriminating Szmielew group. Then
$$\omega>k=\lambda_p+\kappa_{p,n}+\kappa_{p,n+1}+\dots.$$
Hence $\lambda_p<\omega$ and so,
by Proposition~\ref{discr}, $\lambda_p=0$.
Then $$0<\kappa_{p,n}+\kappa_{p,n+1}+\dots<\omega,$$ and so
$I_p$ is finite. Then we have $n<l_p$, and
$\kappa_{p,l_p-1}<\omega$.
In this case $A$ is not discriminating,
by Proposition~\ref{discr}. A contradiction.
\end{proof}

\begin{lemma}
If \textup{(c)} then
$\psi$ holds in some discriminating $p$-Szmielew  group.
\end{lemma}

\begin{proof}
We have $A=\oplus_qA(q)$,
where $A(q)$ is a $q$-Szmielew group.
Put
$$A'(p):=A(p)\oplus\Zpi^{(\omega)}.$$
By Proposition~\ref{discr}, $A'(p)$
is a discriminating $p$-Szmielew group.
Moreover,  $A'(p)\models\psi$.
Indeed, for any sentence $\theta$ of one of the forms
$$\Phi_k(p,n),\ \Phi^k(p,n),\
\Theta^k(p,n),\
\Gamma_k(p,n),\ \Gamma^k(p,n),\
\Delta^k(p,n)$$
if $A\models\theta$ then
$A'(p)\models\theta$.
\end{proof}

It remains to consider case (d). We will need

\begin{lemma}
For any $n\ge k$ the sentence
$\Gamma_l(p,k)$ is effectively equivalent
in abelian groups to
a positive Boolean combination
of sentences of the forms
$\Gamma_i(p,n)$ and $\Phi_j(p,s)$,
where $k\le s<n$ and $0\le i,j\le l$.
\end{lemma}

\begin{proof}
It suffices to show that
in abelian groups
$\Gamma_l(p,k)$ is equivalent to
$$
\Gamma'_l(p,k):=
\bigvee_{i=0}^l(\Gamma_{l-i}(p,k+1)\wedge\Phi_i(p,k)).
$$
A Szmielew group $A$ satisfies
$\Gamma_l(p,k)$ if and only if
$$\mu_p+\kappa_{p,k}+\kappa_{p,k+1}+\dots=l;$$
the latter holds if and only if, for some $i\in\{0,1,\dots,l\},$
$$
\mu_p+\kappa_{p,k+1}+\kappa_{p,k+2}+\dots=l-i\quad
\text{and}\quad \kappa_{p,k}=i,$$
which means that
$\Gamma'_l(p,k)$
holds in $A$.
\end{proof}

Let $n<\omega$ be given.
Replace in $\psi$ every conjunct
$\Gamma_l(p,k)$, where $k<n$, with an equivalent
positive Boolean combination
of sentences of the forms
$\Gamma_i(p,n)$ and $\Phi_j(p,s)$.
The resulting formula is equivalent to
a disjunction of $p$-conjunctions in each of which
there is no conjunct
$\Gamma_l(p,k)$ with $k<n$.
Therefore it remains to prove the following statement,
which allows to decide whether $\psi$ holds in
some discriminating $p$-Szmielew group, in case (d).

\begin{lemma}
Suppose that $\psi$ has
\begin{enumerate}
\item[(a)]
a conjunct
$\Delta_1(p,n)$ or $\Theta_0(p,n)$;
\item[(b)]
no conjuncts
$\Delta_k(p,m)$ with $k\ne 1$ and
$\Theta_k(p,m)$ with $k>0$;
\item[(c)]
no
conjuncts $\Gamma_l(p,s)$ with $s<n$.
\end{enumerate}
Then  the following are equivalent:
\begin{enumerate}
\item
$\psi$ fails in any discriminating $p$-Szmielew group;
\item
there exist  $m$ with $m<n$  and
 $i>0$
such that
\begin{enumerate}
\item[(i)]
$\Phi_i(p,m)$ is a conjunct of $\psi$,
\item[(ii)]
for every $k$ with $m<k<n$
there is $j$ such that
$\Phi_j(p,k)$
is a conjunct of $\psi$.
\end{enumerate}
\end{enumerate}
\end{lemma}

\begin{proof}
First we show that (b) implies that
$\psi$ holds in some $p$-Szmielew group.
If $\Delta_1(p,n)$ is in $\psi$ then
$p^nA=0$; therefore
$A$ is a direct sum of cyclic $p$-groups and hence
a $p$\nobreakdash-Szmielew group.
Suppose $\Delta_1(p,n)$ is not in $\psi$.
Let $A=\oplus_qA(q)$,
where  each $A(q)$ is a $q$-Szmielew group.
Since $\psi$ is a $p$\nobreakdash-conjunction without
conjuncts of the form
$\Delta_k(p,n)$, the $p$-Szmielew group $A(p)\oplus\bQ$
satisfies $\psi$.

So we may assume that $A$ is a $p$-Szmielew group.
By (a), $$\lambda_p=\kappa_{p,n}=\kappa_{p,n+1}\dots=0.$$
Indeed,
if  $\Delta_1(p,n)$ is in $\psi$ then $p^nA=0$;
if  $\Theta_0(p,n)$ is in $\psi$ then
$$0=\lambda_p+\kappa_{p,n}+\kappa_{p,n+1}+\dots.$$
In particular, the set $I_p$ is finite.

Suppose (2).
Due to (i), we have $\kappa_{p,m}=i>0$, and therefore $m<l_p\le n$.
Let $m<k<n$.
By (ii) $\psi$ has a conjunct
$\Phi_j(p,k)$; then
$\kappa_{p,k}=j$.
So $\kappa_{p,k}<\omega$ for all $k$ with
$m\le k<n$.
In particular,
$\kappa_{p,\,l_p-1}<\omega$.
By Proposition~\ref{discr}, in this case $A$ cannot be discriminating, and
(1) follows.

Assuming that (2) is not true, we show that (1) is not true, too.

If $I_p=\emptyset$ then $A$ itself
is discriminating, by Proposition~\ref{discr}.

Suppose $I_p\ne\emptyset$.
First we show that there is $k<n$ such that
$\kappa_{p,r}=0$ for $r>k$, and
for every $j$ the sentence $\Phi_j(p,k)$
is not a conjunct of $\psi$.
Let $m=l_p-1$ and
$i=\kappa_{p,m}$. Then $m<n$ and $i>0$.
If (i) fails, put $k:=m$.
If (i) holds then (ii) fails, and therefore
there is $k$ with $m<k<n$ such that
for every $j$ the sentence $\Phi_j(p,k)$
is not a conjunct of $\psi$.

By Proposition~\ref{discr},
the $p$-Szmielew group
$A\oplus\Zpks^{(\omega)}$
is discriminating.
Moreover,
$$A\oplus\Zpks^{(\omega)}\models\psi.$$
Indeed, by (c) and the choice of $k$,
a conjunct $\theta$ of $\psi$
can have only the forms
$$\Phi_j(p,r),\
\Theta_0(p,n),\
\Gamma_j(p,s),\
\Delta_1(p,n),$$
where $r\ne k$ and $s\ge n$, or the forms
$$\Phi^j(p,t),\
\Theta^j(p,t),\
\Gamma^j(p,t),\
\Delta^j(p,t).$$
Therefore $A\models\theta$ implies
$A\oplus\Zpks^{(\omega)}\models\theta$,
for all such $\theta$.
Here we use that $s\ge n>k$ when consider $\theta$
of the forms $\Theta_0(p,n)$ and $\Gamma_j(p,s)$.
\end{proof}

The proof of Theorem~\ref{decidable} is completed.
\end{proof}

\section{Open questions}

\begin{proposition}\label{nilp-undecidable}
The theory of square-like nilpotent groups is undecidable.
\end{proposition}

\begin{proof}
In fact, even
the universal theory of square-like nilpotent groups is undecidable.
Indeed, it coincides with the universal theory
of nilpotent groups because any nilpotent group $G$ embeds in
the discriminating nilpotent group $G^\omega$.
As any finitely generated nilpotent group
is residually finite, the universal theory
of nilpotent groups coincides with
the universal theory of finite nilpotent groups.
The latter is undecidable~\cite{Ha}.
\end{proof}

\begin{question}
Is the theory of square-like $2$-step nilpotent groups undecidable?
\end{question}

Note that the \emph{universal}
theory of square-like 2\nobreakdash-step nilpotent groups is decidable.
Indeed, as above, it coincides with
the universal theory of 2\nobreakdash-step nilpotent groups and
with the universal theory
of finite 2-step nilpotent groups.
Obviously, the universal theory of 2-step nilpotent groups
is computably enumerable, and
the universal theory of finite 2-step nilpotent groups
is co-computably-enumerable; so the result follows.

Thus, undecidability
of the theory of square-like $2$-step nilpotent groups
cannot be shown like in the proof of
Proposition~\ref{nilp-undecidable}.
In~\cite[Theorem~5.1]{B2} we proved undecidability
of the theory of square-like groups by constructing
a discriminating group which interprets the ring of integers.

\begin{question}
Is there a discriminating $2$-step nilpotent group
which interprets the ring of integers?

\end{question}

Existence of such a group would imply undecidability
of the theory of square-like $2$-step nilpotent groups.

\end{document}